# Hankel determinants of Schröder-like numbers

*Johann Cigler*

Fakultät für Mathematik
Universität Wien
A-1090 Wien, Nordbergstraße 15

johann.cigler@univie.ac.at

## Abstract

After a short survey about Schröder numbers and some generalizations which I call Schröder-like numbers I study some $q-$analogues which have simple Hankel determinants.

## 1. Schröder numbers

In the first paragraphs I state some results about Schröder numbers and Hankel determinants which are either well known or simple modifications of well-known results.
The (large) Schröder numbers $R_n$ can be defined by the generating function

$$F(z)=\sum_{n\geq 0}R_n z^n=\frac{1-z-\sqrt{1-6z+z^2}}{2z}. \tag{1.1}$$

The first terms of this sequence are (cf. [4] A006318)

$$\left(R_n\right)_{n\geq 0}=\left(1,2,6,22,90,\cdots\right). \tag{1.2}$$

The generating function satisfies

$$F(z)=1+zF(z)+zF(z)^2. \tag{1.3}$$

A closely related sequence (cf. [4], A001003) is the sequence $(1,1,3,11,45,197,\cdots)$ of little Schröder numbers $S_n$, defined by $S_0=1$ and $S_n=\frac{R_n}{2}$ for $n>0.$
Their generating function

$$f(z)=\frac{1+F(z)}{2} \tag{1.4}$$

satisfies

$$f(z)=1-zf(z)+2zf(z)^2. \tag{1.5}$$

Historical remarks about these numbers can be found in [7].

It is well known that the Hankel determinants of these numbers are (cf. [4], A006318, A001003)

$$\det\left(R_{i+j}\right)_{i,j=0}^{n-1} = 2^{\binom{n}{2}}, \tag{1.6}$$

$$\det\left(R_{i+j+1}\right)_{i,j=0}^{n-1} = 2^{\binom{n+1}{2}}, \tag{1.7}$$

$$\det\left(R_{i+j+2}\right)_{i,j=0}^{n-1} = 2^{\binom{n+1}{2}}\left(2^{n+1}-1\right), \tag{1.8}$$

$$\det\left(S_{i+j}\right)_{i,j=0}^{n-1} = \det\left(S_{i+j+1}\right)_{i,j=0}^{n-1} = 2^{\binom{n}{2}} \tag{1.9}$$

and

$$\det\left(S_{i+j+2}\right)_{i,j=0}^{n-1} = 2^{\binom{n}{2}}\left(2^{n+1}-1\right). \tag{1.10}$$

There are many ways to prove such results. My favorite method uses orthogonal polynomials. In the sequel I consider only sequences of numbers $(a(n))$ with the property that $a(0)=1$ and $\det(a(i+j))_{i,j=0}^{n-1} \neq 0$ for all $n \geq 1$.
In this situation the polynomials

$$p(n,x) = \frac{1}{\det(a(i+j))_{i,j=0}^{n-1}} \det\begin{pmatrix} a(0) & a(1) & \cdots & a(n-1) & 1 \\ a(1) & a(2) & \cdots & a(n) & x \\ a(2) & a(3) & \cdots & a(n+1) & x^2 \\ \vdots & & & & \vdots \\ a(n) & a(n+1) & \cdots & a(2n-1) & x^n \end{pmatrix} \tag{1.11}$$

are well defined and orthogonal with respect to the linear functional $F$ defined by $F(x^n) = a(n)$.

By Favard's theorem they satisfy a recurrence of the form

$$p(n,x) = (x - s(n-1))p(n-1,x) - t(n-2)p(n-2,x). \tag{1.12}$$

Let $a(n,k)$ be the uniquely determined coefficients in the expansion

$$x^n = \sum_{k=0}^{n} a(n,k)p(k,x). \tag{1.13}$$

They satisfy

$$\begin{aligned} &a(0,k) = [k=0] \\ &a(n,0) = s(0)a(n-1,0) + t(0)a(n-1,1) \\ &a(n,k) = a(n-1,k-1) + s(k)a(n-1,k) + t(k)a(n-1,k+1). \end{aligned} \tag{1.14}$$

Obviously

$$a(n,0)=F(x^n)=a(n). \tag{1.15}$$

We say that the sequences $\left(s(n)\right)$ and $\left(t(n)\right)$ are associated with the sequence $(a(n))$. They contain the same information as the sequence $(a(n))$.

For the large Schröder numbers the associated sequences are (see (2.10) )

$$\begin{aligned} &s(0)=2, s(n)=3 \\ &t(n)=2. \end{aligned} \tag{1.16}$$

This gives the Schröder triangle ( [4] A133367 )

2 1

6 5 1

22 23 8 1

90 107 49 11 1

For the little Schröder numbers the generating function (1.5) implies (see (2.11) )

$$s(0)=1, s(n)=3 \tag{1.17}$$

and

$$t(n)=2. \tag{1.18}$$

This gives the triangle

1 1

3 4 1

11 17 7 1

45 76 40 10 1

The numbers $a(n,k)$ have a well known combinatorial interpretation.
Consider lattice paths, so-called Motzkin paths, with upward steps $(n,k)\to(n+1,k+1)$, horizontal steps $(n,k)\to(n+1,k)$ and downward steps $(n,k+1)\to(n+1,k)$ which start in $(0,0)$ and never fall under the $x-$axis. To each path we associate a weight in the following way: Each upward step has weight 1, each horizontal step at height $k$ has weight $s(k)$ and each downward step which ends on height $k$ has weight $t(k)$. The weight of a path is the product of the weights of its steps. Then $a(n,k)$ is the sum of the weights of all paths from $(0,0)$ to $(n,k)$.
Let now $f_j(z)$ be the generating function of all paths which start and end at height $j$ and never fall under this height. Then $f_j(z)$ is the generating function $f_j(z)=\sum_{n\geq 0} a_j(n,0)z^n$, where $a_j(n,k)$ is given by (1.14) with $s_j(n)=s(n+j)$ and $t_j(n)=t(n+j)$.

Then a simple combinatorial argument gives

$$f(z)=f_0(z)=1+s(0)zf(z)+t(0)z^2 f(z)f_1(z). \tag{1.19}$$

Let

$$d(n,k)=\det\left(a(i+j+k)\right)_{i,j=0}^{n-1}. \tag{1.20}$$

be the Hankel determinant of order $k$ of $\left(a(n)\right)$.

Then (cf. e.g. [6])

$$d(n,0)=\prod_{i=1}^{n-1}\prod_{k=0}^{i-1}t(k)=t(0)^{n-1}t(1)^{n-2}\cdots t(n-3)^2 t(n-2) \tag{1.21}$$

and

$$d(n,1)=(-1)^n p(n,0)d(n,0). \tag{1.22}$$

Sometimes the sequence $\left(d(n,0)\right)$ is called the Hankel transform of the sequence $\left(a(n)\right)$. Under the stated requirements it would be better to call $\left(d(n,0)\right)_{n\geq 1}\times\left(d(n,1)\right)_{n\geq 1}$ the Hankel transform of $\left(a(n)\right)$ since the sequence $\left(a(n)\right)$ is uniquely determined by the sequences of Hankel determinants $\left(d(n,0)\right)$ and $\left(d(n,1)\right)$.

It is often convenient to consider a sequence $(a(n))$ which satisfies $a(2n)=c(n)$ and $a(2n+1)=0$. In this case $s(n)=0$ for all $n$. The corresponding lattice paths have no horizontal steps, i.e. they are so called Dyck paths.

Let

$$\left(a_0(n)\right)=\left(a(2n,0)\right)=\left(c(n)\right) \tag{1.23}$$

Then the associated sequences are

$$s_0(0)=t(0), s_0(n)=t(2n-1)+t(2n), t_0(n)=t(2n)t(2n+1) \tag{1.24}$$

For

$$\left(a_1(n)\right)=\left(a(2n+1,1)\right) \tag{1.25}$$

we get

$$s_1(0)=t(0)+t(1), s_1(n)=t(2n)+t(2n+1), t_1(n)=t(2n+1)t(2n+2). \tag{1.26}$$

As an example consider the sequence $(a(n))$ with $a(2n) = R_n$ and $a(2n+1) = 0$. In this case we get $s(n) = 0, t(2n) = 2, t(2n+1) = 1$ (cf. (2.9)).

The corresponding triangle is

0 1

2 0 1

0 3 0 1

6 0 5 0 1

0 11 0 6 0 1

22 0 23 0 8 0 1

Since by (1.14) $a_1(n) = S_{n+1}$ are again little Schröder numbers, (1.26) gives a triangle for the little Schröder numbers $s_{n+1}$. This is [4], A110440.

3 1

11 6 1

45 31 9 1

197 156 60 12 1

The generating function for the sequence $\left(S_{n+1}\right)$ is $g(z) = \dfrac{f(z)-1}{z}$ and satisfies $g(z) = 1 + 3zg(z) + 2z^2 g(z)^2$.

## 2. Schröder-like numbers

We first want to study Schröder-like numbers $A(n, x, y)$ defined by the generating function $F(z) = \sum_{n\geq 0} A(n, x, y) z^n$ which satisfies

$$F(z) = 1 + xzF(z) + yzF(z)^2. \tag{2.1}$$

A useful combinatorial interpretation of $A(n, x, y)$ can be obtained in the following way: Consider lattice paths with upward and downward steps of length 1 and horizontal steps of length 2. If each upward step has weight 1, each downward step has weight $xz$ and each horizontal step has weight $yz$, then the weight of the set of all non-negative paths from $(0,0)$ to $(2n,0)$ is $A(n, x, y)z^n$.

This implies that

$$A(n, x, y) = \sum_{k=0}^{n} \binom{n+k}{2k} C_k x^{n-k} y^k = \sum_{k=0}^{n} \binom{2n-k}{k} C_{n-k} x^k y^{n-k}, \tag{2.2}$$

where $C_k = \dfrac{1}{k+1}\dbinom{2k}{k}$ is a Catalan number.

We first choose a $k-$subset $1 \le a_1 < a_2 < \cdots a_k < 2n$ of $\{1, \cdots, 2n\}$ which contains no successive elements. This is equivalent with $1 \le a_1 < a_2 - 1 < \cdots a_k - k + 1 \le 2n - k.$ Therefore there are $\binom{2n-k}{k}$ possibilities to choose such a subset. Each $a_j$ is the starting point of a horizontal step. On the $2n-2k$ remaining points there are $C_{n-k}$ non-negative Dyck paths.

By solving a quadratic equation we get the explicit formula

$$F(z) = \frac{1 - xz - \sqrt{1 - 2(x+2y)z + x^2 z^2}}{2yz}. \tag{2.3}$$

The numbers $A(n,x,y)$ satisfy the simple recurrence relation

$$A(n,x,y) = \frac{(2n-1)(x+2y)A(n-1,x,y) - (n-2)x^2 A(n-2,x,y)}{n+1} \tag{2.4}$$

with initial values

$$A(0,x,y) = 1 \text{ and } A(1,x,y) = x + y.$$

In order to show this we differentiate (2.1) and get
$F'(z) = xF(z) + xzF'(z) + yF(z)^2 + 2yzF(z)F'(z)$
Substituting (2.1) gives
$zF'(z)(1 - xz - 2yzF(z)) = yzF(z)^2 + xzF(z) = F(z) - 1$

Therefore

$$zF'(z) = \frac{F(z)-1}{(1 - xz - 2yzF(z))} = \frac{F(z)-1}{\sqrt{1-2(x+2y)z+x^2z^2}} = \frac{(F(z)-1)\sqrt{1-2(x+2y)z+x^2z^2}}{1-2(x+2y)z+x^2z^2}$$
$$= \frac{(F(z)-1)\left(1 - xz - 2yzF(z)\right)}{1-2(x+2y)z+x^2z^2} = \frac{\left(xz + 1 + (x+2y)zF(z) - F(z)\right)}{1-2(x+2y)z+x^2z^2}$$

i.e.

$$zF'(z)\left(1 - 2(x+2y)z + x^2z^2\right) = \left(xz + 1 + (x+2y)zF(z) - F(z)\right)$$

Comparing coefficients we get
$nA(n,x,y) + 2(x+2y)(n-1)A(n-1,x,y) + x^2(n-2)A(n-2,x,y)$
$= -A(n,x,y) + (x+2y)A(n-1,x,y)$
and thus (2.4).

The simplest cases occur for $x = 0$ where $A(n,0,y) = C_n y^n$ and for $x + 2y = 0,$ where we get $A(2n+2,2,-1) = 0$ and $A(2n+1,2,-1) = (-1)^n C_n.$

The equation $F(z) = 1 + xzF(z) + yzF(z)^2$ implies

$F(z) - yzF(z)^2 = 1 + xzF(z) = 1 - yzF(z) + (x + y)zF(z)$
or
$F(z) = 1 + (x + y)zF(z)f(z)$ with

$$f(z) = \frac{1}{1 - yzF(z)}. \tag{2.5}$$

Therefore (2.1) is equivalent with

$$\begin{aligned} F(z) &= 1 + (x + y)zF(z)f(z), \\ f(z) &= 1 + yzF(z)f(z). \end{aligned} \tag{2.6}$$

This also implies that

$$f(z) = \frac{x + yF(z)}{x + y}. \tag{2.7}$$

It is easily verified that

$$f(z) = 1 - xzf(z) + (x + y)zf(z)^2. \tag{2.8}$$

If we write $f(z) = \sum_{n \geq 0} a(n, x, y)z^n$ then the numbers $a(n, x, y)$ are a generalization of the little Schröder numbers.

Comparing with (1.19) we see that $F(z^2)$ corresponds to

$$\begin{aligned} &s(n) = 0, \\ &t(2n) = x + y, t(2n + 1) = y. \end{aligned} \tag{2.9}$$

Denote the corresponding $a(n,k)$ by $\alpha(n,k)$.

Since $\alpha(2n,0) = A(n, x, y)$ we get for the original sequence
$\left(A(0, x, y), A(1, x, y), A(2, x, y), \cdots\right)$

$$\begin{aligned} &s(0) = x + y, s(n) = x + 2y \\ &t(n) = y(x + y). \end{aligned} \tag{2.10}$$

Observing that $a(n, x, y) = A(n, -x, x + y)$ we get for the sequence $\left(a(n, x, y)\right)$

$$\begin{aligned} &s(0) = y, s(n) = x + 2y \\ &t(n) = y(x + y). \end{aligned} \tag{2.11}$$

The Hankel determinants are

$$D(n,0)=\det\left(A(i+j,x,y)\right)_{i,j=0}^{n-1}=\left(y(x+y)\right)^{\binom{n}{2}} \tag{2.12}$$

and

$$D(n,1)=\det\left(A(i+j+1,x,y)\right)_{i,j=0}^{n-1}=y^{\binom{n}{2}}\left((x+y)\right)^{\binom{n+1}{2}}. \tag{2.13}$$

The first result is obvious. The second follows from (1.22). Let $r(n)=(-1)^n p(n,0)$. Then

$$\begin{aligned} r(n)&=s(n-1)r(n-1)-t(n-2)r(n-2)\\ &=(x+2y)r(n-1)-y(x+y)r(n-2) \end{aligned}$$

with initial values $r(0)=1$ and $r(1)=(x+y)$.
This gives by induction $r(n)=(x+y)^n$.
By changing $x\to -x, y\to x+y$ we get

$$d(n,0)=\det\left(a(i+j,x,y)\right)_{i,j=0}^{n-1}=\left(y(x+y)\right)^{\binom{n}{2}} \tag{2.14}$$

and

$$d(n,1)=\det\left(a(i+j+1,x,y)\right)_{i,j=0}^{n-1}=y^{\binom{n+1}{2}}\left(x+y\right)^{\binom{n}{2}}. \tag{2.15}$$

In order to compute $d(n,2)=\det\left(a(i+j+2,x,y)\right)_{i,j=0}^{n-1}$ we use the condensation formula (see [6] , (4.1))

$$d(n,0)d(n,2)=d(n+1,0)d(n-1,2)+d(n,1)^2. \tag{2.16}$$

Let $d(n,2)=d(n,1)h(n)$. Then (2.16) leads to

$h(n)=(x+y)h(n-1)+y^n$ with $h(0)=1=\dfrac{x+y-y}{x}$ and $h(1)=x+2y=\dfrac{(x+y)^2-y^2}{x}$

which immediately gives $h(n)=\dfrac{(x+y)^{n+1}-y^{n+1}}{x}$.

Therefore

$$d(n,2)=y^{\binom{n+1}{2}}\left(x+y\right)^{\binom{n}{2}}\frac{(x+y)^{n+1}-y^{n+1}}{x}. \tag{2.17}$$

In the same way we get

$$D(n,2)=y^{\binom{n}{2}}\left(x+y\right)^{\binom{n+1}{2}}\frac{(x+y)^{n+1}-y^{n+1}}{x}. \tag{2.18}$$

## 3. q-analogues of Schröder-like numbers

Barcucci et al. [1] have introduced (large) $q-$Schröder numbers by the generating function

$$F(z)=1+zF(z)+qzF(z)F(qz). \tag{3.1}$$

The first terms are
$(1,1+q,(1+q)(1+q+q^2),(1+q)(1+2q+3q^2+3q^3+q^4+q^5),\cdots)$

We want to consider more generally the $q-$Schröder-like numbers $A(n,x,y)$ with generating function

$$F(z)=1+xzF(z)+yzF(z)F(qz). \tag{3.2}$$

For $(x,y)=(0,1)$ this reduces to $F(z)=1+zF(z)F(qz).$ Therefore $A(n,0,1)=C_n(q)$ are the $q-$Catalan numbers of Carlitz.

### Theorem 1

*Let* $F(z)=\sum_{n\geq 0}A(n,x,y)z^n$ *satisfy the identity* $F(z)=1+xzF(z)+yzF(z)F(qz).$

*Then*

$$D(n,0)=\det\left(A(i+j,x,y)\right)_{i,j=0}^{n-1}=q^{\frac{n(n-1)^2}{2}}y^{\binom{n}{2}}(x+y)^{n-1}(x+qy)^{n-2}\cdots(x+q^{n-2}y), \tag{3.3}$$

$$\begin{aligned}D(n,1)&=\det\left(A(i+j+1,x,y)\right)_{i,j=0}^{n-1}=q^{\binom{n}{2}}(x+y)(x+qy)\cdots(x+q^{n-1}y)D(n,0)\\&=q^{\frac{n^2(n-1)}{2}}y^{\binom{n}{2}}(x+y)^n(x+qy)^{n-1}\cdots(x+q^{n-1}y)\end{aligned} \tag{3.4}$$

*and*

$$D(n,2)=q^{\frac{(n-1)n(n+1)}{2}}y^{\binom{n}{2}}\prod_{j=0}^{n-1}(x+q^jy)^{n-j}\frac{\prod_{j=1}^{n+1}(x+q^{j-1}y)-q^{\binom{n+1}{2}}y^{n+1}}{x}. \tag{3.5}$$

**Remark**

A special case of Theorem 1 has been proved by another method in [3].

An analogue of the little Schröder numbers is given by the generating function

$$f(z)=\sum_{n\geq 0}a(n,x,y)z^n=\frac{x+yF(z)}{x+y}. \tag{3.6}$$

It is easily verified that it satisfies the equation

$$f(z)=1-xzf(qz)+(x+y)zf(z)f(qz) \tag{3.7}$$

and that

$$F(z)=1+(x+y)zF(z)f(qz). \tag{3.8}$$

## Theorem 2

*Let* $f(z)=\sum_{n\geq 0}a(n,x,y)z^n$ *satisfy the identity* $f(z)=1-xzf(qz)+(x+y)zf(z)f(qz)$.

*Then*

$$d(n,0)=\det\left(a(i+j,x,y)\right)_{i,j=0}^{n-1}=q^{3\binom{n}{3}}y^{\binom{n}{2}}\prod_{j=1}^{n-1}\left(x+q^j y\right)^{n-j}, \tag{3.9}$$

$$d(n,1)=\det\left(a(i+j+1,x,y)\right)_{i,j=0}^{n-1}=q^{n\binom{n}{2}}y^{\binom{n+1}{2}}(x+qy)^{n-1}\cdots(x+q^{n-1}y), \tag{3.10}$$

*and*

$$d(n,2)=q^{\frac{(n+1)n(n-1)}{2}}y^{\binom{n+1}{2}}(x+qy)^{n-1}\cdots(x+q^{n-1}y)\frac{\prod_{j=1}^{n+1}(x+q^{j-1}y)-q^{\binom{n+1}{2}}y^{n+1}}{x}. \tag{3.11}$$

For the proofs I use an idea which I already have used in [2]:

There is a uniquely determined series $h(z)=h(z,y)=1+\sum_{n\geq 1}h_n z^n$ such that

$$F(z)=F(z,y)=\frac{h(qz,y)}{h(z,y)}. \tag{3.12}$$

From the defining equation for $F(z)$ we get $\frac{h(qz)}{h(z)}=1+xz\frac{h(qz)}{h(z)}+yz\frac{h(qz)}{h(z)}\frac{h(q^2z)}{h(qz)}$ and therefore $h(qz)=h(z)+xzh(qz)+yzh(q^2z)$.

Comparing coefficients we get $(q^n-1)h_n=q^{n-1}\left(x+q^{n-1}y\right)h_{n-1}$.

This implies

$$h(z,y)=\sum_{k\geq 0}q^{\binom{k}{2}}\frac{(x+y)(x+qy)\cdots(x+q^{k-1}y)}{(q-1)(q^2-1)\cdots(q^k-1)}z^k. \tag{3.13}$$

Thus

$(x+y)h(z,qy)=xh(z,y)+yh(qz,y)$,

i.e.

$$\frac{h(z,qy)}{h(z,y)}=\frac{x+yF(z)}{x+y}=f(z) \tag{3.14}$$

and

$$h(qz,y)-h(z,y)=\sum_{k\geq 0} q^{\binom{k}{2}} \frac{(x+y)(x+qy)\cdots(x+q^{k-1}y)}{(q-1)(q^2-1)\cdots(q^k-1)}(q^k-1)z^k$$

$$=(x+y)z\sum_{k\geq 0} q^{\binom{k-1}{2}} \frac{(x+qy)\cdots(x+q^{k-1}y)}{(q-1)(q^2-1)\cdots(q^{k-1}-1)}(qz)^{k-1}=(x+y)zh(qz,qy)$$

i.e.

$$F(z)=1+(x+y)z\frac{h(qz,qy)}{h(z,y)}. \tag{3.15}$$

Now we can prove the following analogue of (2.6)

$$\begin{aligned} F(z,y)&=1+(x+y)zF(z,y)f(qz,y) \\ f(z,y)&=1+yzf(z,y)F(z,qy). \end{aligned} \tag{3.16}$$

For we have

$$F(z,y)f(qz,y)=\frac{h(qz,y)}{h(z,y)}\frac{h(qz,qy)}{h(qz,y)}=\frac{h(qz,qy)}{h(z,y)}=\frac{h(z,qy)}{h(z,y)}\frac{h(qz,qy)}{h(z,qy)}=f(z,y)F(z,qy).$$

Therefore (3.15) implies the first line in (3.16).

From (3.14) we get

$$f(z,y)=\frac{(x+y)+(yF(z,y)-y)}{x+y}=1+y\frac{F(z,y)-1}{x+y}=1+yz\frac{h(qz,qy)}{h(z,y)}=1+yzf(z,y)F(z,qy).$$

From (3.16) we can deduce the associated sequences $s(n)$ and $t(n)$ for $F(z^2)$.
We get

$$s(n)=0, t(2n)=q^n(x+q^n y), t(2n+1)=q^{2n+1}y. \tag{3.17}$$

This again implies the associated sequences for $F(z)$ by (1.24).
They are

$$\begin{aligned} s(0)&=x+y, s(n)=q^n(x+q^{n-1}y(1+q)) \\ t(n)&=q^{3n+1}y(x+q^n y). \end{aligned} \tag{3.18}$$

From this (3.3) follows immediately from (1.21).

In order to show (3.4) let $r(n)=(-1)^n p(n,0)$. This gives
$r(n)=q^{n-1}(x+q^{n-2}(1+q)y)r(n-1)-q^{3n-5}y(x+q^{n-2}y)r(n-2)$
with initial values $r(0)=1$ and $r(1)=x+y$.

It has to be shown that $r(n)=q^{\binom{n}{2}}(x+y)\cdots(x+q^{n-1}y),$

i.e.

$$q^{\binom{n}{2}}(x+y)\cdots(x+q^{n-1}y)=q^{n-1}(x+q^{n-2}(1+q)y)q^{\binom{n-1}{2}}(x+y)\cdots(x+q^{n-2}y)$$
$$-q^{3n-5}y(x+q^{n-2}y)q^{\binom{n-2}{2}}(x+y)\cdots(x+q^{n-3}y).$$

But this is easily verified.

From (3.16) we get for $f(z^2)$

$$s(n)=0,t(2n)=q^{2n}y,t(2n+1)=q^n(x+q^{n+1}y). \tag{3.19}$$

This implies that the associated sequences for $f(z)$ are

$$s(0)=y,s(n)=q^{n-1}(x+q^{n-1}(1+q)y),t(n)=q^{3n}y(x+q^{n+1}y). \tag{3.20}$$

This implies formulas (3.9) and (3.10).

Finally we want to compute the Hankel determinants $D(n,2)$ and $d(n,2)$.

Let $h(n)=\dfrac{D(n,2)}{q^{(n-1)n(n+1)}y^{\binom{n}{2}}\prod_{j=0}^{n-1}(x+q^jy)^{n-j}}.$

Then the condensation formula $D(n,2)D(n,0)=D(n+1,0)D(n-1,2)+D(n,1)^2$ gives

$h(n)=q^nyh(n-1)+\prod_{j=0}^{n-1}(x+q^jy)$ with $h(0)=1.$

It is easily verified that

$$h(n)=\frac{\prod_{j=0}^{n}(x+q^jy)-q^{\binom{n+1}{2}}y^{n+1}}{x}. \tag{3.21}$$

This immediately gives

$$D(n,2)=q^{\frac{(n-1)n(n+1)}{2}}y^{\binom{n}{2}}\prod_{j=0}^{n-1}(x+q^jy)^{n-j}\left(\prod_{j=0}^{n}(x+q^jy)-q^{\binom{n+1}{2}}y^{n+1}\right) \tag{3.22}$$

The computation of $d(n,2)$ follows in the same manner.